\documentclass[10pt]{article} 
\usepackage{amsfonts}
\usepackage{amssymb}
\usepackage{amstext}
\usepackage{amsbsy}

\usepackage{theorem}
\newtheorem{theo}{}
\newtheorem{theor}[theo]{Theorem}
\newtheorem{prop}[theo]{Proposition}
\newtheorem{coro}[theo]{Corollary}

\newtheorem{lemma}[theo]{Lemma}
\newtheorem{example}[theo]{Example}
\newtheorem{remark}[theo]{Remark}

\newenvironment{proof}
	{\par {\bf Proof:}}
 	{\hfill $\square$ \medskip}

\newcommand{\slim}{\mathop{\mbox{\rm s-lim}}}

\def\ds{\displaystyle}
\def\b{\begin{equation}}
\def\e{\end{equation}}
\def\F{{\mathcal F}}
\def\Fc{{\mathfrak F}}
\def\FF{\mathbb F}
\def\G{{\mathcal G}}
\def\CC{{\mathfrak C}}
\def\H{{\mathcal H}}
\def\M{{\mathfrak M}}
\def\NN{{\mathfrak N}}
\def\V{{\mathcal V}}
\def\W{{\mathcal W}}
\def\N{\mathbb N}
\def\Z{\mathbb Z}
\def\R{\mathbb R}
\def\C{\mathbb C}
\def\D{\mathbb D}
\def\I{\mathbb I}
\def\J{\mathbb J}
\def\om{\omega}
\def\<{\langle}
\def\>{\rangle}

\def\bu{\bold u}
\def\bv{\bold v}
\def\bw{\bold w}
\def\bof{\bold f}
\def\bphi{{\boldsymbol{\varphi}}}
\def\bpsi{{\boldsymbol{\psi}}}

\begin{document}

\title{Spectral Models for Orthonormal Wavelets\\ and Multiresolution Analysis of $L^2({\mathbb R})$}

\author{F. G\'omez-Cubillo$^1$ and Z. Suchanecki$^2$}

\maketitle
\begin{center}
{\small
$^1$Dpto de An\'alisis Matem\'atico, Universidad de Valladolid, 
Facultad de Ciencias, 47011 Valladolid, Spain. 
e-mail: {\tt fgcubill@am.uva.es}.\\
$^2$Universit\'e du Luxembourg, Campus Kirchberg, 6 Rue Coudenhove-Kalergi, L-1359 Luxembourg
and Institute of Mathematics and Informatics, University of Opole, Poland.
e-mail: {\tt zdzislaw.suchanecki@uni.lu}
}

\end{center}

\begin{abstract}
Spectral representations of the dilation and translation operators on $L^2(\R)$ are built through appropriate bases. Orthonormal wavelets and multiresolution analysis are then described in terms of rigid operator-valued functions defined on the functional spectral spaces. The approach is useful for computational purposes.
\end{abstract}



%\tableofcontents

%-------------

\section{Introduction}

Wavelets were introduced in the beginning of the 1980s synthesizing ideas originated in engineering (subband coding in signal processing, pyramidal and multipole algorithms in image processing), physics (coherent states, renormalization group), numerical analysis (spline approximation) and pure mathematics (Littlewood-Paley and Calder\'on-Zygmund theories), and provide us with a particularly simple and powerful tool in mathematical analysis with a great variety of applications \cite{D92,HW96,JMR01,LMR97,M98}.

Two unitary operators on $L^2(\R)$ play an important role in wavelet theory: {\it dilation} by $r>0$, $ D_r$, defined by $[ D_r f](x):=r^{1/2}\,f(rx)$, and {\it translation} by $t$, $T_t$, defined by $[T_tf](x):=f(x-t)$.
Wavelet bases are built by application of translations and dilations to an appropriate function.
Since the main stream of development leading up to wavelets departs from Fourier analysis, Fourier transforms and series are involved in most of the techniques in wavelet theory.
A clear reason for this is the fact that the complex exponentials are eigenfunctions of the derivative operator and then give us an spectral representation of the unitary one-parameter group generated by it, i.e., the group of translations $\{T_t:t\in\R\}$ on $L^2(\R)$ \cite{BS}.

Paying attention to the discrete theory of wavelets, the values of $r$ and $t$ are fixed (usually $r=2$ and $t=1$) and one works with the entire powers of the operators $D:=D_2$ and $T:=T_1$. The spectra of $D$ and $T$ coincide with the unit circle $\partial\D$ of the complex plane $\C$ and have constant denumerable multiplicity. Thus, the functional spectral representations of both operators live in direct integrals $L^2(\partial\D,\H)$ of $\H$-valued functions over $\partial\D$, where $\H$ is an auxiliary separable Hilbert space \cite{BS,HP74,VN49}. Direct integrals have been already considered in abstract wavelet theory --see \cite{P04} and references therein--.
Here all these facts are proved in Section \ref{s2}.
Explicit spectral representations of $D$ and $T$ are given in Propositions \ref{pdft} and \ref{ptft}. They are built on the basis of orthonormal bases of $L^2(\R)$ with appropriate structure and are useful for computational purposes. In practice, any orthonormal basis of $L^2(\R)$, after adapting it, can be used to build such spectral representations. To illustrate this the details for the exponential and Haar bases are included in an Appendix at the end of the paper. Moreover, an additional spectral representation for $T$ is defined in Proposition \ref{ptft1} on the basis of the usual Fourier transform. Remark \ref{rm40} and Proposition \ref{p7} evidence that classical results in wavelet theory join our approach under this spectral model.

The main purpose of this work is to characterize orthonormal wavelets and multiresolution analysis (MRA) of $L^2(\R)$ by means of rigid (operator-valued) functions defined on the functional spectral spaces of $D$ and $T$.
A rigid function on $L^2(\partial\D,\H)$ is a function $A$ from $\partial\D$ into the space ${\mathcal L}(\H)$ of bounded operators on $\H$ such that $A(\om)$ is a partial isometry with the same initial subspace for almost every (a.e.) $\om\in\partial\D$. 
 For orthonormal wavelets this is done in Section \ref{s4} (Theorem \ref{th3}) and for MRA in Section \ref{s5} (Theorems \ref{th4} and \ref{th5}).
The key to obtain these characterizations is the close connection between wandering and invariant subspaces of $L^2(\partial\D,\H)$ and rigid and range functions defined on it \cite{HALMOS61,HELSON64,L61,S63}.
The fundamental concepts and results of this theory are collected in Section \ref{swisrrf}.

Section \ref{s5.1} is devoted to get the pair of discrete quadrature mirror filters associated with an MRA from the two-scale relations in this context. Remark \ref{r26} makes clear these filters act just on the generic spectral models given in Proposition \ref{ptft}.

Theorem \ref{th3-cc}, Example \ref{example1}, Corollary \ref{coro19} and Proposition \ref{prop27} highlight the usefulness for computational purposes of the spectral models considered in this work. Finite sets of non-zero coordinates imply local character in ``time" and ``frequency", the meaning of these concepts depending on the bases considered. At the same time, suitable choices of the bases lead to sparse sets of conditions (compare, for example, the matrices $\big(\alpha_{i,n}^{s,j,m}\big)$ for the exponential and Haar bases in Appendix). A detailed analysis of this set of conditions shall be done elsewhere.

%---------------------
\section{Spectral models of dilations and translations}\label{s2}

Let ${\D}$ denote the open unit disc of the complex plane $\C$ and $\partial {\D}$ its boundary:
$$
{\D}:=\{\lambda\in\C:|\lambda|<1\}\,,\qquad \partial {\D}:=\{\omega\in\C:|\omega|=1\}\,.
$$
In $\partial {\D}$ interpret measurability in the sense of Borel and consider the normalized Lebesgue measure $d\omega/(2\pi)$.
Given a separable Hilbert space $\H$, let $L^2(\partial {\D};\H)$ denote the set of all measurable functions $\bv:\partial {\D}\to \H$ such that 
$$\int_{\partial {\D}} ||\bv(\omega)||^2_\H\,\frac{d\omega}{2\pi}<\infty$$ (modulo sets of measure zero); measurability here can be interpreted either strongly or weakly, which amounts to the same due to the separability of $\H$. The functions in $L^2(\partial {\D};\H)$ constitute a Hilbert space with pointwise definition of linear operations and inner product given by 
$$
(\bu,\bv):=\int_{\partial {\D}} \big(\bu(\omega),\bv(\omega)\big)_\H\,\frac{d\omega}{2\pi}\,,\qquad \big(\bu,\bv\in L^2(\partial {\D};\H)\big)\,.
$$
The space $L^2(\partial {\D};\H)$ is a particular case of direct integral of Hilbert spaces. The theory  of direct integrals is originally due to von Neumann \cite{VN49} and is in the basis of the functional spectral models for operators \cite[Chapter 7]{BS}. In the case of a constant field of Hilbert spaces, $\H_\om=\H$ for all $\om\in\partial {\D}$, the direct integral $\int_{\partial {\D}}^\oplus \H_\om\,\frac{d\om}{2\pi}$ is just the space $L^2(\partial {\D};\H)$.  
Elementary properties of vector and operator valued functions can be found in \cite[Chapter III]{HP74}; see also \cite{HALMOS61} and \cite[Lecture VI]{HELSON64}.

Now, let $r>1$ and consider the {\it dilation operator} $ D_r :L^2(\R)\to L^2(\R)$ defined by
$$
[ D_r f](x):=r^{1/2} f(rx)\,,\quad (f\in L^2(\R))\,.
$$
$ D_r $ is unitary on $L^2(\R)$ and then its spectrum is included in $\partial\D$.
In order to give the spectral resolution of $ D_r $, we must find a direct integral of Hilbert spaces
$\int_{\partial {\D}}^\oplus \H_\om\,d\mu(\om)$ and a unitary transformation $\G$ from $L^2(\R)$ onto $\int_{\partial {\D}}^\oplus \H_\om\,d\mu(\om)$ such that
\b\label{di}
[\G  D_r f](\om)=\om\,[\G f](\om)\,,\quad \text{ for a.e }\om\in\partial {\D}\text{ and }f\in L^2(\R)\,.
\e 
In other words, $\G D_r\G^{-1}$ is the operator ``multiplication by $\omega$" in the direct integral:
$$
[\G D_r\G^{-1}\bof](\om)=\om\cdot\bof(\om)\,,\quad \text{ for a.e }\om\in\partial {\D}\text{ and }\bof\in \int_{\partial {\D}}^\oplus \H_\om\,\frac{d\om}{2\pi}\,.
$$ 
For it, let $\{K_{\pm,j}^{(0)}(x)\}_{j\in\J}$ be an orthonormal basis (ONB) of $L^2[\pm 1,\pm r)$, where $\J$ is a denumerable set of indices (usually, $\J=\N$ or $\J=\Z$), so that the set
$$
\{K_{\pm,j}^{(m)}(x):=r^{m/2}K_{\pm,j}^{(0)}(r^mx)\}_{j\in\J}
$$
is an ONB of $L^2[\pm r^{-m},\pm r^{-m+1})$ for every $m\in\Z$.
For each $f\in L^2(\R)$ put
$$
\tilde f_{\pm,j}^{(m)}:=\int_{\pm r^m}^{\pm r^{m+1}} f(x)\, \overline{K_{\pm,j}^{(m)}(x)}\,dx\,,\quad (j\in\J,\,m\in\Z)\,.
$$
Let $l^2(\J)$ denote the Hilbert space of sequences of complex numbers $(c_j)_{j\in\J}$ such that $\sum_{j\in\J} |c_j|^2<\infty$, and let $\big\{u_{s,j}\big\}_{s=\pm,j\in\J}$ be a fixed ONB of $l^2(\J)\oplus l^2(\J)$, where $\oplus$ denotes orthogonal sum.

\begin{prop}\label{pdft}
The operator $\G$ defined by
\b\label{dft}
\begin{array}{rccl}
\G: & L^2(\R) & \longrightarrow & \ds \int_{\partial {\D}}^\oplus [l^2(\J)\oplus l^2(\J)]\,\frac{d\om}{2\pi}=L^2\big(\partial {\D};l^2(\J)\oplus l^2(\J)\big) 
\\[3ex]
& f & \mapsto & \ds \tilde \bof:=\bigoplus_{s=\pm}\bigoplus_{j\in\J} \left[ \sum_{m\in\Z} \om^{m}\,\tilde f_{s,j}^{(m)} \right]\,u_{s,j}\,.  
\end{array}
\e
determines a functional spectral model for the dilation operator $D_r$, i.e., $\G$ is unitary and satisfies (\ref{di}).
\end{prop}

\begin{proof}
$\G$ is unitary since $\{K_{s,j}^{(m)}(x)\}_{s=\pm,j\in\J,m\in\Z}$ is an ONB of $L^2(\R)$ and $||f||=||\tilde\bof||^2=\sum_{s,j,m} |\tilde f_{s,j}^{(m)}|^2$ for every $f\in L^2(\R)$. Moreover, $\G$ satisfies (\ref{di}) because $\widetilde{[ D_r f]}_{s,j}^{(m)}=\tilde f_{s,j}^{(m-1)}$ for all $s=\pm,\,j\in\J,\,m\in\Z$.
\end{proof}

In a similar way, for $t>0$ and the {\it translation operator} $T_t$ on $L^2(\R)$ defined by
$$
[T_t f](x):=f(x-t)\,,\quad (f\in L^2(\R))\,,
$$
given an ONB $\{L_{i}^{(0)}(x)\}_{i\in\I}$ of $L^2[0,t)$, with $\I$ a denumerable set of indices, the translated set
$$
\{L_{i}^{(n)}(x):=L_{i}^{(0)}(x-tn)\}_{i\in\I}
$$
is an ONB of $L^2[tn,t(n+1))$ for every $n\in\Z$. For each $f\in L^2(\R)$ write 
$$
\hat f_{i}^{(n)}:=\int_{tn}^{{t(n+1)}} f(x)\,\overline{L_{i}^{(n)}(x)}\,dx\,,\quad (i\in\I,\,n\in\Z)\,,
$$
and let $\big\{u_{i}\big\}_{i\in\I}$ be a fixed ONB of $l^2(\I)$.

\begin{prop}\label{ptft}
The operator $\F$ given by
\b\label{tft}
\begin{array}{rccl}
\F: & L^2(\R) & \longrightarrow & \ds \int_{\partial {\D}}^\oplus l^2(\I)\,\frac{d\om}{2\pi}=L^2\big(\partial {\D};l^2(\I)\big)
\\[3ex]
& f & \mapsto & \hat \bof:=\ds \bigoplus_{i\in\I} \left[ \sum_{n\in\Z} \om^{n}\,\hat f_{i}^{(n)} \right]\,u_i\,,
\end{array}
\e
determines a functional spectral model for the translation operator $T_t$, i.e., $\F$ is unitary and satisfies
\b\label{dit}
[\F T_t f](\om)=\om\,[\F f](\om)\,,\quad \text{ for a.e }\om\in\partial {\D}\text{ and }f\in L^2(\R)\,.
\e 
\end{prop}

\begin{proof}
The operator $\F$ is unitary since $\{L_{i}^{(n)}(x)\}_{i\in\I,n\in\Z}$ is an ONB of $L^2(\R)$ and $||f||=||\hat\bof||^2=\sum_{i,n} |\hat f_{i}^{(n)}|^2$ for every $f\in L^2(\R)$. $\F$ satisfies (\ref{dit}) because $\widehat{[ T_t f]}_{i}^{(n)}=\tilde f_{i}^{(n-1)}$ for all $i\in\I,\,n\in\Z$.
\end{proof}

Since $\{K_{s,j}^{(m)}(x)\}_{s=\pm,j\in\J,m\in\Z}$ and $\{L_{i}^{(n)}(x)\}_{i\in\I,n\in\Z}$ are ONB of $L^2(\R)$, given $f\in L^2(\R)$, one has (in $L^2$-sense)
\b\label{fdobl}
f=\sum_{s,j,m} \tilde f^{(m)}_{s,j} K_{s,j}^{(m)}\quad \text{ and }\quad f=\sum_{i,n} \hat f^{(n)}_{i} L_i^{(n)}\,.
\e
The change of representation between the two models is governed by the matrix $\big(\alpha_{i,n}^{s,j,m}\big)$, where
$$
\alpha_{i,n}^{s,j,m}:=\int_\R L_i^{(n)}(x)\,\overline{K_{s,j}^{(m)}(x)}\,dx\,,
$$
so that
$$
K_{s,j}^{(m)}=\sum_{i,n} \overline{\alpha_{i,n}^{s,j,m}}\,L_i^{(n)}\,,\quad
L_i^{(n)}=\sum_{s,j,m} \alpha_{i,n}^{s,j,m}\, K_{s,j}^{(m)}
$$ 
and
\b\label{crf}
\tilde f_{s,j}^{(m)}= \sum_{i,n} \alpha_{i,n}^{s,j,m}\, \hat f^{(n)}_{i}\,,\qquad
\hat f^{(n)}_{i}=\sum_{s,j,m} \overline{\alpha_{i,n}^{s,j,m}} \,\tilde f_{s,j}^{(m)}\,.
\e

In what follows we shall focus attention on the values $r=2$ and $t=1$, i.e., the operators
$$
D:=D_2\quad \text{ and }\quad T:=T_1\,.
$$
For these values, the elements of the matrix $\big(\alpha_{i,n}^{s,j,m}\big)$ for the exponential and Haar bases are given in the Appendix.

Fixed  the orthonormal bases $\{K_{\pm,j}^{(m)}(x)\}_{j\in\J,m\in\Z}$ and $\{L_{i}^{(n)}(x)\}_{i\in\I,n\in\Z}$ of $L^2(\R)$ giving rise to the respective spectral models (\ref{dft}) and (\ref{tft}) of $D$ and $T$, for $f\in L^2(\R)$ we shall write
$$
\F f=\hat{\bold f}=\big\{\hat f_{i}^{(n)}\big\}\,,\qquad 
\G f=\tilde{\bold f}=\big\{\tilde f_{s,j}^{(m)}\big\}\,.
$$

\begin{lemma}
Let $f\in L^2(\R)$ and, for $p,q\in\Z$,
$$
g:= D^pT^q\,f\,,\qquad h:=T^q D^p\,f\,.
$$
If $\F f=\big\{\hat f_{i}^{(n)}\big\}$, $\F g=\big\{\hat g_{i}^{(n)}\big\}$, $\F h=\big\{\hat h_{i}^{(n)}\big\}$, $\G f=\big\{\tilde f_{s,j}^{(m)}\big\}$, $\G g=\big\{\tilde g_{s,j}^{(m)}\big\}$ and $\G h=\big\{\tilde h_{s,j}^{(m)}\big\}$, then
\b\label{qtpd}
\hat g_{l}^{(k)}=\sum_{s,j,m} \overline{\alpha_{l,k}^{s,j,m}}\,\sum_{i,n} \alpha_{i,n+q}^{s,j,m-p}\,\hat f_{i}^{(n)}\,,\quad (l\in\I,\,k\in\Z)\,,
\e
\b\label{pdqt}
\hat h_{l}^{(k)}=\sum_{s,j,m} \overline{\alpha_{l,k-q}^{s,j,m}}\,\sum_{i,n} \alpha_{i,n}^{s,j,m-p}\,\hat f_{i}^{(n)}\,,\quad (l\in\I,\,k\in\Z)\,,
\e
\b\label{qtpd-1}
\tilde g_{r,l}^{(k)}=\sum_{i,n} \alpha_{i,n}^{r,l,k-p}\,\sum_{s,j,m} \overline{\alpha_{i,n-q}^{s,j,m}}\,\tilde f_{s,j}^{(m)}\,,\quad (r=\pm,\,l\in\J,\,k\in\Z)\,,
\e
\b\label{pdqt-1}
\tilde h_{r,l}^{(k)}=\sum_{i,n} \alpha_{i,n}^{r,l,k}\,\sum_{s,j,m} \overline{\alpha_{i,n-q}^{s,j,m+p}}\,\tilde f_{s,j}^{(m)}\,,\quad (r=\pm,\,l\in\J,\,k\in\Z)\,.
\e
\end{lemma}

\begin{proof}
The result is a straightforward consequence of the change of representation formulas (\ref{crf}) and the fact that $\G D\G^{-1}$ is the operator ``multiplication by $\om$" on $\G L^2(\R)$ and $\F D\F^{-1}$ is the operator ``multiplication by $\om$" on $\F L^2(\R)$.
\end{proof}

Finally, let $\hat f$ denote the usual {\it Fourier transform} of a function $f\in L^2(\R)$:
\b\label{FT}
\hat f(y):=\int_\R f(x)\,e^{-2\pi ixy}\,dx\,,\quad (y\in\R)\,.
\e

\begin{prop}\label{ptft1}
Let $\big\{u_{k}\big\}_{k\in\Z}$ be a fixed ONB of $l^2(\Z)$. The operator $\F_*$ defined by  
\b\label{tft1}
\begin{array}{rccl}
\F_*: & L^2(\R) & \longrightarrow & \ds L^2\big(\partial \D;l^2(\Z)\big)
\\[2ex]
& f & \mapsto & \ds \hat \bof_*:=\bigoplus_{k\in\Z} \,\hat f_k(\om)\,u_k\,,
\end{array}
\e
where, if $\om=e^{2\pi i\theta}$, 
$$
\hat f_k(\om)=\hat f_k(e^{2\pi i\theta}):=\overline{\hat f(\theta+k)}\,,\quad \text{ for a.e. }\theta\in[0,1)\text{ and }k\in\Z\,,
$$
determines a functional model for the translation operator $T$.
\end{prop}

\begin{proof}
The operator $\F_*$ is unitary since, for $f\in L^2(\R)$,
$$
||f||^2=||\hat f||^2=\sum_{k\in\Z} \int_{0}^{1} |\hat f(\theta+k)|^2\,d\theta=\sum_{k\in\Z} \int_{\partial\D} |\hat f_k(\om)|^2\,\frac{d\om}{2\pi}=||\hat\bof_*||^2\,.
$$
Moreover, if $\om=e^{2\pi i\theta}$, then 
$$
[\widehat{Tf}]_k(\om)=\overline{[\widehat{Tf}](\theta+k)}=\overline{e^{-2\pi i(\theta+k)}\hat f(\theta+k)}=\om\cdot\hat f_k(\om)\,,
$$
so that $[\F_*Tf](\om)=\om\cdot[\F_*f](\om)$, for a.e. $\om\in\partial\D$ and $f\in L^2(\R)$.
\end{proof}

%--------------------------

\section{Wandering and invariant subspaces. Rigid and range functions.}\label{swisrrf}

Let $\H$ be a separable Hilbert space and denote by ${\mathcal L}(\H)$ the space of bounded linear operators on $\H$.
Consider the functional space $L^2(\partial {\D};\H)$ defined in Section \ref{s2} and the subspace $\CC$ of $L^2(\partial {\D};\H)$ consisting of all constant functions, i.e., the functions $\bv:\partial\D\to\H$ such that there exists a vector $v\in\H$ with $\bv(\om)=v$ for a.e. $\om\in\partial {\D}$.
From now on the operator ``multiplication by $\om$" on $L^2(\partial {\D};\H)$ shall be denoted by 
$M$, that is, 
\b\label{fshift}
[M\bv](\om):=\om\cdot \bv(\om)\,,\quad (\bv\in L^2(\partial {\D};\H),\, \om\in\partial {\D})\,.
\e
It is easy to verify that $M$ is unitary and $M^{-1}(=M^*)$ is defined by $[M^*\bv](\om)=\om^*\cdot \bv(\om)$.

A subspace of a Hilbert space is called a {\it wandering subspace} for an operator $U$ if it is orthogonal  to all its images under the (positive) powers of $U$. For isometries, wandering subspaces behave better than usual: if $U$ is an isometry and if $\M$ is a wandering subspace for $U$, then $U^m\M\perp U^n\M$ whenever $m$ and $n$ are distinct non-negative integers.
If $U$ is unitary, even more is true: in that case $U^m\M\perp U^n\M$ whenever $m$ and $n$ are any two distinct integers (possibly negative).

A weakly measurable\footnote{That $A$ is {\it weakly measurable} means that the scalar product $(A(\omega)h,g)_\H$ is a Borel measurable scalar function on $\partial\D$ for each $h,g\in\H$.} operator-valued function 
$$
A:\partial {\D}\to{\mathcal L}(\H):\om\mapsto A(\om)
$$
is called a {\it rigid operator function} if $A(\omega)$ is for a.e. $\omega\in \partial {\D}$ a partial isometry\footnote{An operator $B\in{\mathcal L}(\H)$ is a {\it partial isometry} if there is a (closed) subspace $\M$ of $\H$ such that $||Bu||=||u||$ for $u\in\M$ and $Bv=0$ for $v\in\M^\perp$. In such case $\M$ is called the {\it initial space} of $B$.} on $\H$ with the same initial space.

According to Halmos \cite[Lemma 5]{HALMOS61}, wandering subspaces for $M$ and rigid functions are related as follows:

\begin{lemma}\label{lwr} {\rm [Halmos]}
A subspace $\M$ of $L^2(\partial {\D};\H)$ is a wandering subspace for $M$ if and only if there exists a rigid function $A$ such that $\M=A\,\CC$.
The subspace $\M$ uniquely determines $A$ to within a constant partially isometric factor on the right.
\end{lemma}

On the other hand, a closed subspace $\M$ of $L^2(\partial {\D};\H)$ is called {\it invariant} if $Mv\in\M$ for every $v\in\M$.
$\M$ is called {\it doubly invariant} if $Mv$ and $M^{-1}v$ belong to $\M$ for each $v\in\M$.
$\M$ is called {\it simply invariant} if it is invariant but not doubly invariant.

A {\it range function} $J=J(\omega)$ is a function on $\partial {\D}$ taking values in the family of closed subspaces of $\H$. $J$ is said {\it measurable} if the orthogonal projection $P(\omega)$ on $J(\omega)$ is weakly measurable. Range functions which are equal a.e. on $\partial {\D}$ are identified. For each measurable range function $J$, $\M_J$ denotes the set of all functions $v$ in $L^2(\partial {\D};\H)$ such that $v(\omega)$ lies in $J(\omega)$ a.e.. The correspondence between $J$ and $\M_J$ is one-to-one, under the convention that range functions are identified if they are equal a.e..

The following result is implicitly contained in the work of Lax \cite{L61} and explicitly proved by Srinivasan \cite{S63} --see also \cite[Theorem 8]{HELSON64}--.

\begin{lemma}\label{ldir}{\rm [Lax-Srinivasan]}
\it The doubly invariant subspaces of $L^2(\partial {\D};\H)$ are the subspaces $\M_J$, where $J$ is a measurable range function.
\end{lemma}

In what follows by the {\it range} of a doubly invariant subspace $\M_J$ we mean that range function $J$. This definition can be extended to an arbitrary set of functions: the {\it range} of a set of vector functions is the range of the smallest doubly invariant subspace containing all the functions.

To determine the simply invariant subspaces of $L^2(\partial {\D};\H)$ one needs to introduce the {\it Hardy classes}.
We denote by $H^2({\D};\H)$ the Hardy class of functions
$$
\tilde \bu(\lambda)=\sum_{k=0}^\infty \lambda^k a_k,\qquad (\lambda\in {\D},\,a_k\in\H),
$$
with values in $\H$, holomorphic on ${\D}$, and such that
$\frac{1}{2\pi}\int_{\partial {\D}} ||\tilde \bu(r\omega)||^2_\H\,d\omega$, ($0\leq r<1$),
has a bound independent of $r$ or, equivalently, such that $\sum ||a_k||^2_\H<\infty$.
For each function $\tilde \bu\in H^2({\D};\H)$ the non-tangential limit in strong sense 
$$
\slim_{\lambda\to\omega} \tilde \bu(\lambda)=\sum_{k=0}^\infty \omega^k a_k=: \bu(\omega)
$$
exist for almost all $\omega\in {\partial {\D}}$. The functions $\tilde \bu(\lambda)$ and $ \bu(\omega)$ determine each other (they are connected by Poisson formula), so that we can identify $H^2({\D};\H)$ with a subspace of $L^2(\partial {\D};\H)$, say $H^+(\partial {\D};\H)$, thus providing $H^2({\D};\H)$ with the Hilbert space structure of $H^+(\partial {\D};\H)$ and embedding it in $L^2(\partial {\D};\H)$ as a subspace.\footnote{ 
Fixing some orthonormal basis $\{u_1,u_2,\cdots\}$ for $\H$, each element $\bof$ of $L^2(\partial {\D};\H)$ may be expressed in terms of its relative coordinate functions $f_j$,
$$
\bof(\omega)=\sum_{j}f_j(\omega)\cdot u_j\,,
$$
where the functions $f_j$ are in $L^2(\partial {\D}):=L^2(\partial {\D};\C)$ and $||\bof||^2=\sum_{j} ||f_j||^2$. A function $\bof$ in $L^2(\partial {\D};\H)$ belongs to $H^+(\partial {\D};\H)$ iff its coordinate functions $f_j$ are all in the Hardy space $H^+(\partial {\D}):=H^+(\partial {\D};\C)$.
}
The complementary space  $H^-(\partial {\D};\H):=L^2(\partial {\D};\H)\backslash H^+(\partial {\D};\H)$ is associated with the conjugate Hardy class $H^2(\C\backslash \overline{{\D}};\H)$ in a similar way.

The general form of the simply invariant subspaces of $L^2(\partial {\D};\H)$ is given by Helson \cite[Theorem 9]{HELSON64} on the basis of the work of Halmos \cite[Theorems 3 and 4]{HALMOS61}. 

\begin{lemma}\label{lsis}{\rm [Halmos-Helson]}
Each simply invariant subspace $\M$ of $L^2(\partial {\D};\H)$ has the form
\begin{equation}\label{H64.6.31}
\M=A\, H^+(\partial {\D};\H)\oplus \M_K\,,
\end{equation}
where $K$ is a measurable range function and $A$ is a rigid function with range $J$ orthogonal to $K$ almost everywhere. The doubly invariant part $\M_K$ of $\M$ is just $\cap_{n\in\Z} M^n\M$.
The subspace $\M$ uniquely determines the rigid operator function $A$ to within a constant partially isometric factor on the right.
\end{lemma}

%------------------

\section{Orthonormal wavelets on $\R$}\label{s4}

An {\it orthonormal wavelet on $\R$} is a function $\psi\in L^2(\R)$ such that $\{\psi_{j,k}:j,k\in\Z\}$ is an orthonormal basis of $L^2(\R)$, where
$$
\psi_{j,k}(x):=[D^j T^k \psi](x)=2^{j/2}\psi(2^j x-k)\,,\quad (j,k\in\Z)\,.
$$ 

\begin{remark}
\label{rm40}
\rm
It is a well known result --see, for example, \cite[Lemma 2.2.4]{LMR97}-- that for $\psi\in L^2(\R)$, $\big\{\psi(\cdot-k):k\in\Z\big\}$ is an orthonormal system if and only if 
$\sum_{k\in\Z} |\hat \psi(\theta+k)|^2= 1$ for a.e. $\theta\in\R$, where $\hat \psi$ is the Fourier transform of $\psi$ defined in (\ref{FT}).
This result fits in with the concept of rigid function with one-dimensional initial subspace on the spectral model for $T$ given by (\ref{tft1}). 
Indeed, $\big\{\psi(\cdot-k):k\in\Z\big\}$ is an orthonormal system in $L^2(\R)$ if and only if the closed subspace spanned by $\hat\bpsi_*:=\F_*\psi$ is a wandering subspace for the unitary operator $M=\F_*T\F_*^{-1}$ defined in $\F_* L^2(\R)$ by (\ref{fshift}). By Lemma \ref{lwr}, this is equivalent to the existence of a rigid function $\hat A$ on $\F_* L^2(\R)$ with one-dimensional initial subspace and such that $\hat\bpsi_*(\om)=\hat A(\om)\, u$, with $u$ a normalized vector belonging to the initial subspace of $\hat A(\om)$, the same one-dimensional subspace of $l^2(\Z)$ for a.e. $\om\in\partial\D$. Since $\hat A$ is a partial isometry for a.e. $\om=e^{2\pi i\theta}\in\partial {\D}$,
$$
\begin{array}{rl}
1&\ds =||u||_{l^2(\Z)}^2=||\hat A(\om) u||_{l^2(\Z)}^2= ||\hat\bpsi_*(\om)||_{l^2(\Z)}^2=\\[2ex]
&\ds =\sum_{k\in\Z} |\hat \psi_k(e^{2\pi i\theta})|^2= \sum_{k\in\Z} |\hat \psi(\theta+k)|^2\,, \quad \text{for a.e. }\theta\in\R\,.
\end{array}
$$
\end{remark}

From now on in this Section we shall pay attention to the respective spectral models for $D$ and $T$ given in Propositions \ref{pdft} and \ref{ptft}.
$\tilde \CC$ shall denote the subspace of $\G L^2(\R)=L^2\big(\partial {\D};l^2(\J)\oplus l^2(\J)\big)$ of all constant functions and $\hat \CC$ the subspace of constant functions of $\F L^2(\R)=L^2\big(\partial {\D};l^2(\I)\big)$, that is, 
$$
\tilde \CC:= \G\big(L^2(-r,-1]\oplus L^2[1,r)\big)\,,\quad \hat \CC:= \F\big(L^2[0,1)\big)\,.
$$
The expression $\<\cdot\>$ shall denote the closed subspace expanded by the argument ``$\cdot$":
$$
\<\cdot\>:=\overline{\text{span}\{\cdot\}}\,.
$$
We shall use the {\it discrete Dirac delta function} $\delta_k$ defined by
$$
\delta_k:=\left\{\begin{array}{ll}1,&\text{ if }k=0,\\ 0,&\text{ if }k\in\Z\backslash\{0\}.\end{array}\right.
$$

The next result is straightforward.

\begin{lemma}\label{loww}
Let $\psi$ be an orthonormal wavelet of $L^2(\R)$.
Then:
\begin{itemize}
\item[(a)]
$\<\psi\>$ is a wandering subspace for $D$ and $T$ on $L^2(\R)$. Equivalently, $\G\<\psi\>$ is a wandering subspace for $\G D\G^{-1}$ on $\G L^2(\R)$ and $\F\<\psi\>$ is a wandering subspace for $\F T\F^{-1}$ on $\F L^2(\R)$
\item[(b)]
$\<T^k\psi:k\in\Z\>$ is a wandering subspace for $D$ on $L^2(\R)$. Equivalently, $\G\< T^k\psi:k\in\Z\>$ is a wandering subspace for $\G D\G^{-1}$ on $\G L^2(\R)$.
\item[(c)]
$\G\<T^k\psi:k\in\Z\>$ is a doubly invariant subspace of $\F L^2(\R)$.
\end{itemize}
\end{lemma}

Now, the results of Section \ref{swisrrf} allow us to characterize any orthonormal wavelet of $L^2(\R)$ in terms of rigid functions on $\G L^2(\R)$ and $\F L^2(\R)$.

\begin{theor}\label{th3}
If $\psi$ is an orthonormal wavelet of $L^2(\R)$, then there exist rigid functions $\hat A$ on $\F L^2(\R)$ and $\tilde B$, $\tilde C$ on $\G L^2(\R)$ such that:
\begin{itemize}
\item[(i)] the initial subspaces of $\hat A$, $\tilde B$ are one-dimensional, and $\tilde C(\om)$ is unitary for a.e. $\om\in\partial\D$; 
\item[(ii)] $\F^{-1}\<\om^k\,\hat A\,\hat\CC:k\in\Z\>=\G^{-1}\<\tilde C\,\tilde\CC\>$;
\item[(iii)] $\<\psi\>=\F^{-1}\<\hat A\,\hat\CC\>=\G^{-1}\<\tilde B\,\tilde\CC\>$.
\end{itemize}
The wavelet $\psi$ uniquely determines $\hat A$, $\tilde B$ and $\tilde C$ to within constant partially isometric factors on the right.

Conversely, assume that there exist rigid functions $\hat A$ on $\F L^2(\R)$ and $\tilde B$, $\tilde C$ on $\G L^2(\R)$ satisfying (i), (ii) and
(iii') $\F^{-1}\<\hat A\,\hat\CC\>=\G^{-1}\<\tilde B\,\tilde\CC\>$.
Then, the triplet $(\hat A,\tilde B,\tilde C)$ has associated a unique orthonormal wavelet $\psi\in L^2(\R)$ given by 
\b\label{th3-1}
\psi:=\F^{-1}[\hat A\,\hat\bu]=\G^{-1}[\tilde B\,\tilde\bv]\,,
\e
where $\hat\bu(\om)=u$ and $\tilde\bv(\om)=v$ for a.e. $\om\in\partial\D$, being $u$ and $v$ normalized vectors in the initial spaces of $\hat A(\om)$ and $\tilde B(\om)$, respectively. 

Additionally, there exists a measurable range function $\hat J=\hat J(\om)$ on $\F L^2(\R)$ such that $\hat J(\om)$ is one-dimensional for a.e. $\om\in\partial\D$ and
\b\label{th3-2}
\F^{-1}\M_{\hat J}=\F^{-1}\<\om^k\,\hat A\,\hat\CC:k\in\Z\>=\G^{-1}\<\tilde C\,\tilde\CC\>\,.
\e
\end{theor}

\begin{proof}
Let $\psi$ be an orthonormal wavelet of $L^2(\R)$.
According to Lemma \ref{loww}.(a) and Lemma \ref{lwr}, there exist rigid functions $\hat A$ on $\F L^2(\R)$ and $\tilde B$ on $\G L^2(\R)$ satisfying (iii). $\hat A$ and $\tilde B$ have one-dimensional initial subspaces because $\<\psi\>$ is one-dimensional.
Since $\F T\F^{-1}$ is the operator ``multiplication by $\om$" on $\F L^2(\R)$, Lemma \ref{loww}.(b) and Lemma \ref{lwr} ensure that there exists  a rigid function $\tilde C$ on $\G L^2(\R)$ verifying (ii). Moreover, $\tilde C(\om)$ is unitary for a.e. $\om\in\partial\D$ because $\G D\G^{-1}$ is the operator ``multiplication by $\om$" on $\G L^2(\R)$ and $\{D^j T^k\psi:j,k\in\Z\}$ is a basis of $L^2(\R)$, so that
$$
\<\om^j\,\G T^k\psi:j,k\in\Z\>=\G\<D^j T^k\psi:j,k\in\Z\>=\G L^2(\R)\,,
$$
i.e., $\<\tilde C\,\tilde\CC\>$ has full range --see \cite[Lecture VII]{HELSON64}--.
Obviously the function $\psi$ uniquely determines the subspaces $\<\psi\>$ and $\<T^k\psi:k\in\Z\>$. Thus, by Lemma \ref{lwr}, $\psi$ uniquely determines the rigid functions $\hat A$, $\tilde B$ and $\tilde C$ to within constant partially isometric factors on the right.

Conversely, let $(\hat A,\tilde B,\tilde C)$ be a triplet of rigid functions satisfying (i), (ii) and (iii'), and define $\psi\in L^2(\R)$ by (\ref{th3-1}). Then, by Lemma \ref{lwr}, $\{D^j T^k\psi:j,k\in\Z\}$ is an orthogonal system because $\hat A$, $\tilde B$ and $\tilde C$ are rigid functions. One has $||D^j T^k\psi||=1$ for every $j,k\in\Z$, since $D,\,T$ are unitary operators, $\hat A(\om)$ is a partial isometry for a.e. $\om\in\partial\D$ and $u\in l^2(\I)$ is a normalized vector. The completeness of the orthonormal system $\{D^j T^k\psi:j,k\in\Z\}$ follows from the unitarity of $\tilde C(\om)$ for a.e. $\om\in\partial\D$ as before.

Finally, Lemma \ref{loww}.(c) and Lemma \ref{ldir} imply that there exists a range function $\hat J$ satisfying (\ref{th3-2}) and such that $\hat J(\om)=\<[\F\psi](\om)\>=\<\hat A(\om)\,u\>=\<\tilde B(\om)\,v\>$ for a.e. $\om\in\partial\D$.
\end{proof}

\begin{remark}\rm
The rigid functions $\tilde B$ and $\tilde C$ of Theorem \ref{th3} are related as follows:
Let $\NN_1$ and $\NN_2$ be the respective (constant) initial spaces of  $\tilde B$ and $\tilde C$. There exists a rigid function $\tilde G$ on $\G L^2(\R)=L^2\big(\partial {\D};l^2(\J)\oplus l^2(\J)\big)$, with initial space $\NN_1$ and the range of $\tilde G(\om)$ included in $\NN_2$ for a.e. $\om\in\partial\D$, such that $\tilde G\tilde\CC\subseteq H^+\big(\partial {\D};l^2(\J)\oplus l^2(\J)\big)$ and $\tilde B=\tilde C\,\tilde G$. 
See \cite[Lemma 6]{HALMOS61} for details. 
\end{remark}

For practical purposes it is advisable to deal with coordinates.

\begin{theor}\label{th3-cc}
A function $\psi\in L^2(\R)$ with $\G\psi=\tilde\bpsi=\big\{\tilde \psi_{s,j}^{(m)}\big\}$
is an orthonormal wavelet if and only if the following two conditions are satisfied:
\begin{itemize}
\item[(i)]
Orthonormality: For every $p,q\in\Z$,
\b\label{oncw}
\sum_{s,j,m} \tilde\psi_{s,j}^{(m)}\,\sum_{i,n} \overline{\alpha_{i,n}^{s,j,m-p}}\,\sum_{r,k,l}\alpha_{i,n-q}^{r,k,l}\,\overline{\tilde\psi_{r,k}^{(l)}}=\delta_p\,\delta_q\,.
\e
\item[(ii)]
Completeness: For every finite set $\FF$ of indices $(s,j)$, with $s=\pm$ and $j\in\J$, the matrix\footnote{Each file of this matrix corresponds with fixed values of the pair of indices $(m,q)\in\Z\times\Z$, whereas each of its columns corresponds with fixed values of the pair $(s,j)\in\FF$, so that the matrix has an infinite number of files and the number of its columns is the cardinal of $\FF$.} 
$$
\left(\sum_{i,n} \overline{\alpha_{i,n}^{s,j,m}}\,\sum_{r,k,l}\alpha_{i,n-q}^{r,k,l}\,\overline{\tilde\psi_{r,k}^{(l)}}\right)_{\begin{array}{l}\scriptstyle{(m,q)\in\Z\times\Z}\\[-1ex] \scriptstyle{(s,j)\in\FF}\end{array}} 
$$
has maximal range, i.e., the cardinal of $\FF$. 
\end{itemize}
\end{theor}

\begin{proof}
By definition, a function $\psi\in L^2(\R)$ is an orthonormal wavelet if and only the system $\{\psi_{p,q}=D^p T^q \psi:p,q\in\Z\}$ is an orthonormal basis of $L^2(\R)$.
Since $D$ and $T$ are unitary operators, the orthonormality of this system is equivalent to
\b\label{woc}
(\psi_{p,q},\psi)=\delta_p\,\delta_q\,,\quad (p,q\in\Z)\,.
\e
Putting $\G\psi_{p,q}=\tilde \bpsi_{p,q}=\big\{[\tilde \psi_{p,q}]_{s,j}^{(m)}\big\}$, condition (\ref{woc}) can be written as
\b\label{woc-1}
\int_{\partial\D} \om^p\,\big(\tilde\bpsi_{0,q}(\om),\tilde\bpsi(\om)\big)\,\frac{d\om}{2\pi}=\delta_p\,\delta_q\,,\quad (p,q\in\Z)\,.
\e
If $p\neq0$, this means that $\big(\tilde\bpsi_{0,q}(\om),\tilde\bpsi(\om)\big)$, which is a complex-valued integrable function, is such that all its Fourier coefficients are equal to zero. If $p=0$, this means that $\big(\tilde\bpsi_{0,q}(\om),\tilde\bpsi(\om)\big)$ has Fourier coefficients equal to zero, except the constant term, which is $1$. That is, (\ref{woc-1}) is equivalent to
\b\label{woc-2}
\big(\tilde\bpsi_{0,q}(\om),\tilde\bpsi(\om)\big)=\delta_q\,,\quad \text{ for a.e. }\om\in\partial\D\,.
\e
From (\ref{dft}),
$$
\begin{array}{rl}
\ds\big(\tilde\bpsi_{0,q}(\om),\tilde\bpsi(\om)\big)
&\ds =\sum_{s,j}\Big[{\sum_{m} \om^{m}\,\tilde \psi_{s,j}^{(m)}}\Big]
\Big[\overline{\sum_{l} \om^{l}\,[\tilde \psi_{0,q}]_{s,j}^{(l)}}\Big]
\\
&\ds =\sum_p \om^p\,\Big[\sum_{s,j,m} \tilde \psi_{s,j}^{(m)}\overline{[\tilde \psi_{0,q}]_{s,j}^{(m-p)}}\Big]\,.
\end{array}
$$
Thus, in terms of coordinates, (\ref{woc-2}) is written as
$$
\sum_{s,j,m} \tilde \psi_{s,j}^{(m)}\overline{[\tilde \psi_{0,q}]_{s,j}^{(m-p)}}=\delta_p\,\delta_q\,,\quad (p,q\in\Z)\,.
$$
Since, by (\ref{qtpd-1}),
$[\tilde \psi_{0,q}]_{s,j}^{(m-p)}=\sum_{i,n} \alpha_{i,n}^{s,j,m-p}\,\sum_{r,k,l} \overline{\alpha_{i,n-q}^{r,k,l}}\,\tilde \psi_{r,k}^{(l)}$,
one gets the condition (i) of the result.

Now, for $f\in L^2(\R)$ let us write $\G f=\tilde\bof=\big\{\tilde f_{s,j}^{(m)}\big\}$. The subspace $\Fc$ consisting of all functions $f\in L^2(\R)$ with a finite number of non-zero coordinates $f_{s,j}^{(m)}$ is dense in $L^2(\R)$. Thus, the completeness of the orthonormal system $\{\psi_{p,q}:p,q\in\Z\}$ in $L^2(\R)$ is equivalent to the condition: if $f\in\Fc$ and $(f,\psi_{p,q})=0$ for all $p,q\in\Z$, then $f=0$.
From (\ref{fdobl}) and (\ref{qtpd-1}),
$$
(f,\psi_{p,q})=\sum_{s,j,m} \tilde f_{s,j}^{(m)}\,\overline{[\tilde \psi_{p,q}]_{s,j}^{(m)}}=
\sum_{s,j,m}\tilde f_{s,j}^{(m)}\,\sum_{i,n} \overline{\alpha_{i,n}^{s,j,m-p}}\,\sum_{r,k,l}\alpha_{i,n-q}^{r,k,l}\,\overline{\tilde\psi_{r,k}^{(l)}}\,.
$$
Since $m$ and $p$ are independent, one arrives at the condition (ii) of the result.
\end{proof}

\begin{remark}\rm
Rearranging sums in (\ref{oncw}) and using (\ref{crf}) one gets
$$
\sum_{i,n}\, [\hat\psi_{p,0}]_i^{(n)}\,\overline{[\hat\psi_{0,q}]_i^{(n)}}=(\psi_{0,q},\psi_{p,0})=\delta_p\,\delta_q\,,\quad (p,q\in\Z)\,,
$$
where $\F\psi_{p,q}=\big\{[\hat\psi_{p,q}]_i^{(n)}\big\}$. This gives a direct proof for the orthonormality condition (i) of Theorem \ref{th3-cc}.
We include the former proof as an exercise of style in this context.
\end{remark}

\begin{example}\label{example1}\rm
For both the exponential bases of Appendix \ref{appA1} and the Haar bases of Appendix \ref{appA2} one has the following result:
Let $\psi$ be a function of $L^2(\R)$ with compact support included in the interval $[1,2]$ (any other interval of length one could be considered). Then, the only possible non-zero coordinates of $\psi$ are $\tilde\psi_{+,j}^{(0)}, j\in\J$, and $\psi$ is an  orthonormal wavelet if and only if
$$
\sum_{j} \overline{\tilde\psi_{+,j}^{(0)}}\,\sum_{k} {\alpha_{k,1+q}^{+,j,-p}}\,{\tilde\psi_{+,k}^{(0)}}=\delta_p\,\delta_q\,,\quad (p,q\in\Z)\,,
$$
and for every finite set $\FF$ of indices $(s,j)$, with $s=\pm$ and $j\in\J$, the matrix
$$
\left(\sum_{k} {\alpha_{k,1+q}^{s,j,m}}\,{\tilde\psi_{+,k}^{(0)}}\right)_{\begin{array}{l}\scriptstyle{(m,q)\in\Z\times\Z}\\[-1ex] \scriptstyle{(s,j)\in\FF}\end{array}} 
$$
has range the cardinal of $\FF$. It is easier to deal with these conditions when considering the Haar bases of Appendix \ref{appA2}, because the matrix $\big(\alpha_{i,n}^{s,j,m}\big)$ is really sparse for them. The solutions of this set of conditions shall be studied elsewhere.
\end{example}

%------------------

\section{Multiresolution analysis on $\R$}\label{s5}

A {\it multiresolution analysis} ({\it MRA}) of $L^2(\R)$ is a sequence $\{\V_n\}_{n\in\Z}$ of closed subspaces of $L^2(\R)$ such that 
\begin{itemize}
\item[(i)]
 $\{0\}\subset \cdots\subset \V_{-1}\subset \V_0\subset\V_1\subset \cdots\subset L^2(\R)$,
\item[(ii)]
 $\cap_{n\in\Z} \V_n=\{0\}$,
\item[(iii)]
$\overline{\cup_{n\in\Z} \V_n}=L^2(\R)$,
\item[(iv)]
 $f\in\V_0 \Leftrightarrow D^n f\in\V_n$ for all $n\in\Z$,
\item[(v)]
 there is a function $\varphi\in L^2(\R)$, the {\it scaling function}, whose integer translates 
$\{T^k\varphi:k\in\Z\}$ form an orthonormal basis of $\V_0$.\footnote{Sometimes (v) is replaced by the weaker condition:
(v')  there is a function $\varphi\in L^2(\R)$ such that the set $\{T^k\varphi:k\in\Z\}$ is a {\it Riesz basis} of $\V_0$, i.e. its linear span is dense in $\V_0$ and there exist positive constants $A$ and $B$ such that 
$A\sum_{k\in\Z} c_k^2 \leq ||\sum_{k\in\Z} c_k\varphi(\cdot-k) ||^2_{L^2(\R)}\leq B\sum_{k\in\Z} c_k^2$,
for each $\{c_k\}\subset\R$ such that $\sum_{k\in\Z} c_k^2<\infty$.
In fact, one can orthogonalize the basis and (v) and (v') are equivalent. See \cite[Theorem 2.2.5]{LMR97}, \cite[Section 2.1]{HW96}) or \cite[Theorem 7.1]{M98} for details.
}
\end{itemize}

\begin{remark}\rm
The properties used to define an MRA are not independent:
Conditions (i), (iv) and (v) [or (v')] imply (ii) --see, for example, \cite[Theorem 2.1.6]{HW96}--.
Obviously, properties (iv) and (v) imply that the set
$$
\big\{\varphi_{n,k}:=D^nT^k\varphi:k\in\Z\big\}
$$
is an orthonormal basis of $\V_n$ for each $n\in\Z$.
\end{remark}

%----

With the modest prerequisite that $\varphi\in L^1(\R)$  and taking into consideration the spectral model given by $\F_*$ in (\ref{tft1}) one can simply construct MRAs:

\begin{prop}\label{p7}
Let $\varphi\in L^2(\R)\cap L^1(\R)$ and such that:
\begin{itemize}
\item[(i)]
$\sum_{k\in\Z} |\hat \varphi_k(\om)|^2=1$ for a.e. $\om\in\partial {\D}$;
\item[(ii)]
$\hat \varphi_k(1)=\delta_k$.
\end{itemize}
Then $\varphi$ is a scaling function for an MRA $\big\{\V_n\big\}_{n\in\Z}$ of $L^2(\R)$.
\end{prop}

\begin{proof} 
According to Remark \ref{rm40} condition (i) is equivalent to the fact that $\big\{\varphi(\cdot-k):k\in\Z\big\}$ is an orthonormal system.
Moreover, since $\hat \varphi_k(e^{2\pi i\theta})=\hat \varphi(\theta+k)$ for $\theta\in[0,1)$ and $k\in\Z$, in terms of the usual Fourier transform $\hat\varphi$, condition (ii) means that $\hat\varphi(0)=1$ and $\hat \varphi(k)=0$ for $k\in\Z\backslash\{0\}$. It is well known --see, for example, \cite[Theorem 2.2.7]{LMR97}-- that under these conditions the spaces $\V_n:=\<\varphi_{n,k}:k\in\Z\big\>$, $n\in\Z$, give an MRA of $L^2(\R)$.
\end{proof}

Scaling functions and MRAs of $L^2(\R)$ can also be characterized in terms of rigid functions on the spectral models $\G L^2(\R)$ and $\F L^2(\R)$ given in Section \ref{swisrrf}.
The key to do it is the following simple result:

\begin{lemma}\label{lsfw}
Let $\varphi$ be a scaling function for an MRA $\big\{\V_n\big\}_{n\in\Z}$ of $L^2(\R)$. Then:
\begin{itemize}
\item[(a)]
$\<\varphi\>$ is a wandering subspace for $T$ on $L^2(\R)$. Equivalently, $\F\<\varphi\>$ is a wandering subspace for $\F T\F^{-1}$ on $\F L^2(\R)$
\item[(b)]
The orthogonal complement of $\G\V_0=\G\<T^k\varphi:k\in\Z\>$ is a (simply) invariant subspace on $\G L^2(\R)$ without doubly invariant part.
\item[(c)]
$\F\V_0=\F\<T^k\psi:k\in\Z\>$ is a doubly invariant subspace of $\F L^2(\R)$.
\end{itemize}
\end{lemma}

\begin{theor}\label{th4}
Let $\varphi$ be a scaling function for an MRA $\big\{\V_n\big\}_{n\in\Z}$ of $L^2(\R)$. Then there exist rigid functions $\hat R$ on $\F L^2(\R)$ and $\tilde S$ on $\G L^2(\R)$ such that:
\begin{itemize}
\item[(i)] the initial subspace of $\hat R$ is one-dimensional and $\tilde S(\om)$ is unitary for a.e. $\om\in\partial\D$; 
\item[(ii)] $\V_0=\F^{-1}\<\om^k\,\hat R\,\hat\CC:k\in\Z\>=\G^{-1}\<\tilde S\,H^-(\partial {\D}; l^2(\J)\oplus l^2(\J))\>$;
\item[(iii)] $\<\varphi\>=\F^{-1}\<\hat R\,\hat\CC\>$.
\end{itemize}
The scaling function $\varphi$ uniquely determines $\hat R$ and $\tilde S$ to within constant partially isometric factors on the right.

Conversely, assume that there exist rigid functions $\hat R$ on $\F L^2(\R)$ and $\tilde S$ on $\G L^2(\R)$ satisfying (i) and
\begin{itemize}
\item[(ii')] $\F^{-1}\<\om^k\,\hat R\,\hat\CC:k\in\Z\>=\G^{-1}\<\tilde S\,H^-(\partial {\D}; l^2(\J)\oplus l^2(\J))\>$.
\end{itemize}
Then, the pair $(\hat R,\tilde S)$ has associated a unique scaling function $\varphi$ for an MRA $\big\{\V_n\big\}_{n\in\Z}$ of $L^2(\R)$ given by 
\b\label{th4-1}
\varphi:=\F^{-1}[\hat R\,\hat\bw]\,,
\e
where $\hat\bw$ is a (constant) normalized vector in the initial space of $\hat R$. In such case, (ii) and (iii) are satisfied. 

Additionally, there exists a measurable range function $\hat K=\hat K(\om)$ on $\F L^2(\R)$ such that $\hat K(\om)$ is one-dimensional for a.e. $\om\in\partial\D$ and
\b\label{th4-2}
\F^{-1}\M_{\hat K}=\F^{-1}\<\om^k\,\hat A\,\hat\CC:k\in\Z\>=\V_0\,.
\e
\end{theor}

\begin{proof}
Let $\varphi$ be a scaling function for an MRA $\big\{\V_n\big\}_{n\in\Z}$ of $L^2(\R)$.
According to Lemma \ref{lsfw}.(a) and Lemma \ref{lwr}, there exists a rigid function $\hat R$ on $\F L^2(\R)$ satisfying (iii). $\hat R$ has one-dimensional initial subspace because $\<\varphi\>$ is one-dimensional.
Recall that $\F T\F^{-1}$ is the operator ``multiplication by $\om$" on $\F L^2(\R)$ and $\G D\G^{-1}$ is the operator ``multiplication by $\om$" on $\G L^2(\R)$.  Then, Lemma \ref{lsfw}.(b) and Lemma \ref{lsis} ensure that there exists  a rigid function $\tilde S$ on $\G L^2(\R)$ verifying (ii). Here the doubly invariant subspace $\M_K$ does not appear because $\cap_{n\in\Z} \V_n=\{0\}$. Moreover, $\tilde S(\om)$ is unitary for a.e. $\om\in\partial\D$ because  and $\{D^j T^k\varphi:j,k\in\Z\}$ is a basis of $L^2(\R)$, so that
$\<\tilde S\,\tilde\CC\>$ has full range --see \cite[Lecture VII]{HELSON64}--.
Since the function $\varphi$ uniquely determines the subspaces $\<\varphi\>$ and $\<T^k\varphi:k\in\Z\>$, by Lemmas \ref{lwr} and \ref{lsis}, $\hat R$ and $\tilde S$ are determined by $\varphi$ to within constant partially isometric factors on the right.

Conversely, let $(\hat R,\tilde S)$ be a pair of rigid functions satisfying (i) and (ii'), and define $\varphi\in L^2(\R)$ by (\ref{th4-1}). Then, by Lemma \ref{lwr}, $\{T^k\varphi:k\in\Z\}$ is an orthogonal system because $\hat R$ is a rigid function. One has $||T^k\psi||=1$ for every $k\in\Z$ since $T$ is unitary, $\hat R(\om)$ is a partial isometry for a.e. $\om\in\partial\D$ and $w\in l^2(\I)$ is a normalized vector. 
That $\overline{\cup_{n\in\Z} \V_n}=L^2(\R)$ follows from the unitarity of $\tilde S(\om)$ for a.e. $\om\in\partial\D$ as before. That $\G\V_0$ has no doubly invariant part implies $\cap_{n\in\Z} \V_n=\{0\}$. 

Finally, Lemma \ref{lsfw}.(c) and Lemma \ref{ldir} imply that there exists a range function $\hat K$ satisfying (\ref{th4-2}) and such that $\hat K(\om)=\<[\F\varphi](\om)\>=\<\hat R(\om)\,w\>$ for a.e. $\om\in\partial\D$.
\end{proof}

\begin{coro}\label{coro19}
Let $\varphi\in L^2(\R)$ be a scaling function and 
$\F\varphi=\hat\bphi=\big\{\hat \varphi_{i}^{(n)}\big\}$. Then
\b\label{rofb}
\sum_{i,n} \hat \varphi_{i}^{(n)}\,\overline{\hat \varphi_{i}^{(n-k)}}=\delta_k\,.
\e
\end{coro}

\begin{proof}
Condition (\ref{rofb}) is equivalent to $||\hat\bphi(\om)||=1$ for a.e. $\om\in\partial {\D}$. Indeed,
$$
\begin{array}{c}
\ds ||\hat\bphi(\om)||^2=\sum_n \big|\sum_i \hat \varphi_{i}^{(n)}\,\om^n\big|^2=
\sum_n \big(\sum_i \hat \varphi_{i}^{(n)}\,\om^n\big)\big(\sum_j \overline{\hat \varphi_{i}^{(m)}}\,\om^{-m}\big)=
\\
=\ds \sum_{i,k} \om^k\, \big(\sum_i \hat \varphi_{i}^{(n)}\,\overline{\hat \varphi_{i}^{(n-k)}}\big)=
\sum_k \om^k\,\big(\sum_{i,n} \hat \varphi_{i}^{(n)}\,\overline{\hat \varphi_{i}^{(n-k)}}\big)
\end{array}
$$
and the last expression is equal to $1$ for a.e. $\om\in\partial {\D}$ if and only if (\ref{rofb}) is satisfied.
\end{proof}

Given a MRA, denote by $\W_n$ the orthogonal complement of $\V_n$ in $\V_{n+1}$, i.e. $\V_{n+1}=\W_n\oplus\V_n$. From (ii),
\b\label{vw}
\V_{n+1}=\V_n\oplus\W_n=\bigoplus_{j=-\infty}^n \W_j\,,\quad (n\in\Z)\,.
\e
From (iii),
\b\label{cw}
L^2(\R)=\bigoplus_{n\in\Z} \W_n\,.
\e

Now, suppose that $\psi$ is an orthonormal wavelet. For $j\in\Z$ let $\W_j$ be the closure in $L^2(\R)$ of the span of $\{\psi_{j,k}=D^jT^k\psi:k\in\Z\}$. Since $\psi$ is an orthonormal wavelet, (\ref{cw}) is satisfied. Moreover, the sequence of subspaces $\{\V_n:n\in\Z\}$ defined by (\ref{vw}) satisfies properties (i)--(iv).
Thus,  $\{\V_n:n\in\Z\}$ will generate an MRA if there exists a function $\varphi\in L^2(\R)$ such that the system $\{T^k\varphi:k\in\Z\}$ is an orthonormal basis for $\V_0$. 
In this case we say that the wavelet $\psi$ is {\it associated with and MRA} or, more simply, that $\psi$ is an {\it MRA wavelet}.

When $\psi$ is an MRA wavelet of $L^2(\R)$ with scaling function $\varphi$ the rigid functions $\tilde C$ of Theorem \ref{th3} and $\tilde S$ of Theorem \ref{th4} coincide. Therefore:

\begin{theor}\label{th5}
Let $\psi$ be an MRA wavelet of $L^2(\R)$ with scaling function $\varphi$.
Then there exist rigid functions $\hat A$, $\hat R$ on $\F L^2(\R)$ and $\tilde B$, $\tilde S$ on $\G L^2(\R)$ such that:
\begin{itemize}
\item[(i)] the initial subspaces of $\hat A$, $\tilde B$, $\hat R$ are one-dimensional, and $\tilde S(\om)$ is unitary for a.e. $\om\in\partial\D$; 
\item[(ii)] one has
$$
\begin{array}{l}
\W_0=\F^{-1}\<\om^k\,\hat A\,\hat\CC:k\in\Z\>=\G^{-1}\<\tilde S\,\tilde\CC\>\,,
\\[1ex]
\V_0=\F^{-1}\<\om^k\,\hat R\,\hat\CC:k\in\Z\>=\G^{-1}\<\tilde S\,H^-(\partial {\D}; l^2(\J)\oplus l^2(\J))\>\,;
\end{array}
$$
\item[(iii)] $\<\psi\>=\F^{-1}\<\hat A\,\hat\CC\>=\G^{-1}\<\tilde B\,\tilde\CC\>$ and $\<\varphi\>=\F^{-1}\<\hat R\,\hat\CC\>$.
\end{itemize}
The MRA wavelet $\psi$ and scaling function $\varphi$ uniquely determine $(\hat A,\,\tilde B,\,\hat R,\,\tilde S)$ to within constant partially isometric factors on the right.

Conversely, assume that there exist rigid functions $\hat A$, $\hat R$ on $\F L^2(\R)$ and $\tilde B$, $\tilde S$ on $\G L^2(\R)$ satisfying (i) and
$$
\begin{array}{c}
\F^{-1}\<\om^k\,\hat A\,\hat\CC:k\in\Z\>=\G^{-1}\<\tilde S\,\tilde\CC\>\,,
\\[1ex]
\F^{-1}\<\om^k\,\hat R\,\hat\CC:k\in\Z\>=\G^{-1}\<\tilde S\,H^-(\partial {\D}; l^2(\J)\oplus l^2(\J))\>\,,
\\[1ex]
\F^{-1}\<\hat A\,\hat\CC\>=\G^{-1}\<\tilde B\,\tilde\CC\>\,.
\end{array}
$$
Then the set $(\hat A,\,\tilde B,\,\hat R,\,\tilde S)$ has associated a unique MRA wavelet $\psi\in L^2(\R)$ with scaling function $\varphi$ given by 
\b\label{th5-1}
\psi:=\F^{-1}[\hat A\,\hat\bu]=\G^{-1}[\tilde B\,\tilde\bv]\,,\quad \varphi:=\F^{-1}[\hat R\,\hat\bw]\,,
\e
where $\hat\bu$, $\hat\bv$ and $\hat\bw$ are (constant) normalized vectors in the initial spaces of $\hat A$, $\tilde B$ and $\hat R$, respectively. In such case, (ii) and (iii) are satisfied. 

Additionally, there exist measurable range functions $\hat J$ and $\hat K$ on $\F L^2(\R)$ such that $\hat J(\om)$, $\hat K(\om)$ are one-dimensional for a.e. $\om\in\partial\D$ and
\b\label{th5-2}
\W_0=\F^{-1}\M_{\hat J}\,,\quad \V_0=\F^{-1}\M_{\hat K}\,.
\e
\end{theor}

%-----------

\subsection{Two-scale relations and mirror filters}\label{s5.1}

Whenever an MRA of $L^2(\R)$ is given there exists an orthonormal wavelet  associated with it. 
Moreover, the wavelet can be constructed explicitly from the MRA. This is carry out thanks to the two-scale relations and a pair of discrete quadrature mirror filters appears. 
As we shall see in this Section, these filters act just on the spectral model $\F L^2(\R)$ given in Proposition \ref{ptft}.

Let $\big\{\V_n\big\}_{n\in\Z}$ be an MRA of $L^2(\R)$ with scaling function $\varphi$. Since $\varphi\in\V_0$, $\V_0\subset \V_1$ and $\big\{\varphi_{1,k}:k\in\Z\big\}$ is an orthonormal basis of $\V_1$, there exists a sequence $\{h_k\}\in l^2(\Z)$ such that
\b\label{tsr}
\varphi(x)=\sum_{k\in\Z} h_k\,\varphi_{1,k}(x)=\sum_{k\in\Z} h_k\,\varphi(2x-k)\,.
\e
Equation (\ref{tsr}) is usually known as the {\it two-scale relation} of the scaling function $\varphi$.
A simple change of variable $y=2^jx$ generalizes this relation to any scale:
\b\label{tsrg}
\varphi_{j,0}=\sum_{k\in\Z} h_k\,\varphi_{j+1,k}\,,\quad (j\in\Z)\,.
\e
In terms of the dilation and translation operators, $D$ and $T$, since $\varphi_{j,k}= D^jT^k \varphi$, relation (\ref{tsrg}) is written as
\b\label{tsro}
 D^j\varphi=\sum_{k\in\Z} h_k\, D^{j+1}T^k\varphi\,,\quad (j\in\Z)\,.
\e
Taking $j=-1$ in (\ref{tsro}) and going to the spectral model given in (\ref{tft}), by (\ref{dit}) one obtains
\b\label{tsrs}
[\F D^{-1}\varphi](\om)=\sum_{k\in\Z} h_k\om^k\, [\F\varphi](\om)=h(\om)\,[\F\varphi](\om)\,,
\e
where the {\it two-scale symbol} for $\varphi$, defined by
$$
h(\om):=\sum_{k\in\Z} h_k\,\om^k\in L^2(\partial {\D})\,,
$$
appears in a natural way.

In a similar way, if $\psi$ is a wavelet associated to the MRA $\big\{\V_n\big\}_{n\in\Z}$, since $\psi_{j,0}\in\W_j$ and $\V_{j+1}=\V_j\oplus\W_j$ for $j\in\Z$, {\it two-scale relations for $\psi$} are also satisfied: there is a sequence $\{g_k\}\in l^2(\Z)$ such that
\b\label{tsrgw}
\psi_{j,0}=\sum_{k\in\Z} g_k\,\varphi_{j+1,k}\,,\quad (j\in\Z)\,.
\e
In particular, for $j=-1$,
\b\label{tsrsw}
[\F D^{-1}\psi](\om)=\sum_{k\in\Z} g_k\om^k\, [\F\varphi](\om)=g(\om)\,[\F\varphi](\om)\,,
\e
where $g$ is  the {\it two-scale symbol for $\psi$}:
$$
g(\om):=\sum_{k\in\Z} g_k\,\om^k\in L^2(\partial {\D})\,.
$$

The following result is easily proved using the spectral model for $T$ of Proposition \ref{ptft}. As usual, the symbol $\perp$ means ``orthogonal to", the symbol $\oplus$ ``orthogonal sum" and the {\it overline} ``adherence of" or "closure of".

\begin{prop}\label{p8}
Let $\big\{\V_n\big\}_{n\in\Z}$ be an MRA of $L^2(\R)$ with scaling function $\varphi$ and wavelet $\psi$.
Let $h$ and $g$ be the two-scale symbols for $\varphi$ and $\psi$, respectively. Then $h$ and $g$ are elements of $L^2(\partial {\D})$ satisfying the following properties:
\begin{itemize}
\item[(a)]
$\om^{2n}h(\om)\perp h(\om)$, for all $n\in\Z\backslash\{0\}$;
\item[(b)]
$\om^{2m}g(\om)\perp g(\om)$, for all $m\in\Z\backslash\{0\}$;
\item[(c)]
$\om^{2n}h(\om)\perp \om^{2m}g(\om)$, for all $n,m\in\Z$;
\item[(d)]
$\overline{\oplus_{n\in\Z}\<\om^{2n}h(\om)\>}\oplus\overline{\oplus_{m\in\Z}\< \om^{2m}g(\om)\>}=L^2(\partial {\D})$.
\end{itemize}
\end{prop}

\begin{proof}
The unitarity of $ D$ and the relation $T^k D= D\,T^{2k}$ for all $k\in\Z$ give 
$$
\begin{array}{rl}
\V_{-1}&\ds = D^{-1}\V_0= D^{-1}\,\overline{\oplus_{k\in\Z}\<T^k \varphi\>}\\
&\ds =\overline{\oplus_{k\in\Z}\< D^{-1}T^k D D^{-1}\varphi\>}=
\overline{\oplus_{k\in\Z}\<T^{2k} D^{-1}\varphi\>}\,.
\end{array}
$$
This relation together with (\ref{tsrs}) and (\ref{dit}) lead to
$$
\F\V_{-1}=\overline{\oplus_{k\in\Z} \<\om^{2k}\,h(\om)\,[\F\varphi](\om)\>}\,.
$$
Thus (the symbol $\ominus$ means orthogonal complement),
\b\label{fw-1.1}
\begin{array}{rl}
\F\W_{-1}&\ds =\F\V_0\ominus\F\V_{-1}\\
&\ds=\overline{\oplus_{k\in\Z}\<\om^k\,[\F\varphi](\om)\>}\ominus \overline{\oplus_{k\in\Z} \<\om^{2k}\,h(\om)\,[\F\varphi](\om)\>}\,.
\end{array}
\e
Now, from the two-scale relation (\ref{tsrgw}) for $\psi$ with $j=-1$,
$$
\begin{array}{rl}
\W_{-1}&\ds = D^{-1}\W_0=\overline{\oplus_{k\in\Z}\< D^{-1}T^k  D D^{-1}\psi\>}\\
&\ds=\overline{\oplus_{k\in\Z}\<T^{2k} D^{-1}\psi\>}=
\overline{\oplus_{k\in\Z}\<T^{2k}\,\sum_{k\in\Z} g_k\,T^k\varphi\>}
\end{array}
$$
and
\b\label{fw-1.2}
\F\W_{-1}=\overline{\oplus_{k\in\Z}\<\om^{2k}\,g(\om)\,[\F\varphi](\om)\>}\,.
\e
Since $\hat B$ is a rigid function and $v$ is a normalized vector of its initial subspace, $||[\F\varphi](\om)||_{l^2(\Z)}=1$ for a.e $\om\in\partial {\D}$ and, therefore, equating the right hand sides of (\ref{fw-1.1}) and (\ref{fw-1.2}), one gets the result.
\end{proof}

A straightforward consequence of Proposition \ref{p8} is:

\begin{coro}\label{l11}
The two-scale symbol $h(\om):=\sum_{k\in\Z} h_k\om^k$ satisfies the following two equivalent properties:
\b\label{mcc0ae}
\sum_{k\in\Z} \overline{h_k}\,h_{k+2n}=\delta_n\,,
\e
\b\label{mcc1ae}
|h(\om)|^2+|h(-\om)|^2=2\,,\quad \text{ for a.e. } \om\in\partial {\D}\,.
\e
\end{coro}

\begin{proof}
$\sum_{k\in\Z} |h_k|^2=1$ is a consequence of the two-scale relation (\ref{tsr}) and the normalization of $\varphi_{j,k}$.
That for $n\in\Z\backslash\{0\}$ one has $\sum_{k\in\Z} \overline{h_k}\,h_{k+2n}=0$ follows from property (a) in Proposition \ref{p8}.
The equivalence of (\ref{mcc0ae}) and (\ref{mcc1ae}) follows from writing out the explicit Fourier series for $|h(\om)|^2+|h(-\om)|^2$ --see \cite[page 137]{D92} for details--.
\end{proof}

Equivalent conditions to (a)--(d) in Proposition \ref{p8} on the functions $h$ and $g$ are given by 
van Eijndhoven and Oonincx in \cite{VEO99}:

\begin{prop}\label{p9}
{\rm [van Eijndhoven and Oonincx]}
Two functions $h$ and $g$ of $L^2(\partial {\D})$ satisfy the conditions (a)--(d) in Proposition \ref{p8} if and only if
\begin{itemize}
\item[(1)]
$g(\om)\,\overline{h(\om)}+g(-\om)\,\overline{h(-\om)}=0$, for a.e. $\om\in\partial {\D}$;
\item[(2)]
The matrix $\pmatrix{h(\om) & g(\om)\cr h(-\om) & g(-\om)}$ is invertible for a.e. $\om\in\partial {\D}$.
\end{itemize}
\end{prop}

\begin{remark}\label{rqmf}\rm
From conditions (1) and (2) in Proposition \ref{p9} it is possible to find $g$ once $h$ is known (or to find $h$ once $g$ is known). For example, it can be verified in a straightforward way that possible choices of $g$ are given by
\b\label{afh}
g(\om)= \om^{2m+1} \overline{h(-\om)}\,,
\e
for any $m\in\Z$, so that $g(\om)=\sum_{k\in\Z} (-1)^{1-k}\,\overline{h_{2m+1-k}}\,\om^k$ or
$$
g_k=(-1)^{1-k}\,\overline{h_{2m+1-k}}\,,\quad (k\in\Z)\,.
$$
By (\ref{mcc1ae}), for the choices of $g$ given in (\ref{afh}) the matrix $\frac{1}{\sqrt{2}}\pmatrix{h(\om) & g(\om)\cr h(-\om) & g(-\om)}$ is unitary for a.e. $\om\in\partial {\D}$.
Thus, conditions (1) and (2) in Proposition \ref{p9} together with (\ref{afh}) lead to a pair $(h,g)$ of transfer functions for two (discrete) {\it quadrature mirror filters}. See \cite[Section 5,6]{D92} and \cite[Chapter 3]{JMR01} for details.
\end{remark}

\begin{remark}\rm
Property (c) in Proposition \ref{p8} together with (\ref{afh}) lead to 
$$
\sum_{k\in\Z} (-1)^{1-k}\,h_k\,h_{1-2n-k}=0\,,\quad (n\in\Z)\,.
$$
\end{remark}

\begin{remark}\label{r26}\rm
{\it The two discrete quadrature mirror filters $\{h_k\}$ and $\{g_k\}$ act on the spectral model $\F L^2(\R)$ given in Proposition \ref{ptft}.}  Indeed, if $\F\varphi=\hat\bphi=\big\{\hat \varphi_{i}^{(n)}\big\}$, $\F\varphi_{-1,0}=\hat\bphi_{-1,0}=\big\{[\hat \varphi_{-1,0}]_{i}^{(n)}\big\}$ and $\F\psi_{-1,0}=\hat\bpsi_{-1,0}=\big\{[\hat \psi_{-1,0}]_{i}^{(n)}\big\}$, since
$$
\F\varphi_{0,k}=\F T^k\varphi=\om^k\,\hat\bphi(\om)=\big\{\hat \varphi_{i}^{(n-k)}\big\}\,,\quad (k\in\Z)\,,
$$
the two-scale relations (\ref{tsrg}) and (\ref{tsrgw}) for $j=-1$ are equivalent to
$$
[\hat \varphi_{-1,0}]_{i}^{(n)}=\sum_k h_k\,\hat \varphi_{i}^{(n-k)}\,,\quad
[\hat \psi_{-1,0}]_{i}^{(n)}=\sum_k g_k\,\hat \varphi_{i}^{(n-k)}\,,\quad (i\in\I,\,n\in\Z)\,.
$$
On the other hand, relations (\ref{tsrs}) and (\ref{tsrsw}) can be written as
$$
\hat\bphi_{-1,0}(\om)=h(\om)\cdot\hat\bphi(\om)\,, \quad \hat\bpsi_{-1,0}(\om)=g(\om)\cdot\hat\bphi(\om)\,,\quad \text{for a.e. }\om\in\partial\D\,.
$$
In general, for $j\in\Z$, filtering by $\{h_k\}$ and $\{g_k\}$ the coordinates in $\F L^2(\R)$ of $\varphi_{j,0}$ results in the coordinates of $\varphi_{j-1,0}$ and $\psi_{j-1,0}$, and using the transfer functions $h$ and $g$ as multipliers on $\hat\bphi_{j,0}$ results in $\hat\bphi_{j-1,0}$ and $\hat\bpsi_{j-1,0}$, respectively.   
\end{remark}

\begin{remark}\rm
Let $\big\{\V_n\big\}_{n\in\Z}$ be an MRA  of $L^2(\R)$ with scaling function $\varphi\in\V_0$ and orthonormal wavelet $\psi\in\W_0=\V_1\ominus\V_0$. Let $H^+:=H^+(\partial\D;\C)$ denote the space of scalar-valued Hardy functions. Jorgensen \cite{J96} considers a spectral transform for $D^{-1}$ of the form
$$
\G_*: L^2(\R)\to L^2(\partial\D;\W_0)\,,
$$
so that $\G_*\V_1=H^+(\partial\D,\W_0)$ and $\G_*\W_0=\tilde\CC$ (the subspace of constant functions), and proves that there exist functions $a$, $b_k$, $b_k^*$, ($k\in\N$), in $H^+$ such that
$$
\tilde\bphi=a\cdot\tilde\bpsi\,,\quad
\tilde\bphi_{0,k}= b_k\cdot\tilde\bpsi\,,\quad 
\tilde\bphi_{0,-k}= b_k^*\cdot\tilde\bpsi_{0,-1}\,,\quad (k\in\N)\,,
$$
where now $\tilde\bof:=\G_* f$ for $f\in L^2(\R)$. Moreover, there is an inner function $c$, i.e., $c\in H^+$ and $|c(\om)|=1$ for a.e. $\om\in\partial\D$, such that
$$
\om\,c(\om)\, b_k(\om)=b_k^*(\om)\,,\quad (\om\in\partial\D,\, k\in\N)\,.
$$
The inner function $c$ serves as a {\it spectral invariant}. It separates distinct examples and is trivial precisely for the Haar wavelet.
\end{remark}

Taking into account the two spectral models $\G L^2(\R)$ and $\F L^2(\R)$ of Propositions \ref{pdft} and \ref{ptft}, and the change of representation matrix $\big(\alpha_{i,n}^{s,j,m}\big)$,
one reaches the following identities:

\begin{prop}\label{prop27}
Let $\big\{\V_n\big\}_{n\in\Z}$ be an MRA of $L^2(\R)$ with scaling function $\varphi$ and $\F\varphi=\hat\bphi=\big\{\hat \varphi_{i}^{(n)}\big\}$.
If $h(\om)=\sum_{k\in\Z} h_k\,\om^k$ is the two-scale symbol for $\varphi$, then, for every $p\in\I$ and $q\in\Z$,
\b\label{tsrb4}
\begin{array}{rl}
\ds \hat \varphi_p^{(q)} & \ds= \sum_{s,j,m} \overline{\alpha_{p,q}^{s,j,m}}\, \sum_{i,n} \alpha_{i,n}^{s,j,m-1}\,\sum_kh_{k}\,\hat \varphi_i^{(n-k)}
\\[1ex]
&\ds = \sum_k \sum_{s,j,m} \overline{\alpha_{p,q-k}^{s,j,m}}\, \sum_{i,n} \big(h_{2k}\,\alpha_{i,n}^{s,j,m-1}+
h_{2k+1}\,\alpha_{i,n+1}^{s,j,m-1}\big)\,\hat \varphi_i^{(n)}\,.
\end{array}
\e
Moreover, the function $\psi\in L^2(\R)$ defined by $\F\psi=\big\{\hat \psi_{i}^{(n)}\big\}$, where, for every $p\in\I$ and $q\in\Z$,
\b\label{tsrb4w}
\begin{array}{rl}
\ds \hat \psi_p^{(q)} & \ds= \sum_{s,j,m} \overline{\alpha_{p,q}^{s,j,m}}\, \sum_{i,n} \alpha_{i,n}^{s,j,m-1}\,\sum_k (-1)^{1-k}h_{1-k}\,\hat \varphi_i^{(n-k)}
\\[1ex]
&\ds = \sum_k \sum_{s,j,m} \overline{\alpha_{p,q-k}^{s,j,m}}\, \sum_{i,n} \big(h_{1-2k}\,\alpha_{i,n}^{s,j,m-1}-
h_{-2k}\,\alpha_{i,n+1}^{s,j,m-1}\big)\,\hat \varphi_i^{(n)}\,,
\end{array}
\e
is an orthonormal wavelet associated with the MRA $\big\{\V_n\big\}_{n\in\Z}$.
\end{prop}

\begin{proof}
The two-scale relation $\varphi_{-1,0}=\sum_k h_k\,\varphi_{0,k}$ on $\F L^2(\R)$ looks like
\b\label{tsrb}
\hat\bphi=\hat D(h\cdot\hat\bphi)=\hat D\Big(\sum_k h_k\, \hat T^k\,\hat\bphi\Big)=\sum_k h_k\, \hat D\,\hat T^k\,\hat\bphi\,,
\e
where $\hat D:=\F D\F^{-1}$ and $\hat T:=\F T\F^{-1}$.
By (\ref{qtpd}), the first line in (\ref{tsrb4}) is the counterpart of (\ref{tsrb}) for the coordinates $\big\{\hat \varphi_{i}^{(n)}\big\}$.
Moreover, since $\hat D\,\hat T^{2k}=\hat T^{k}\,\hat D$ for every $k\in\Z$,
\b\label{tsrb1}
\begin{array}{rl}
\hat\bphi=\sum_k h_k\, \hat D\,\hat T^k\,\hat\bphi & 
=\sum_k \big[h_{2k}\,\hat D\,\hat T^{2k}+h_{2k+1}\,\hat D\,\hat T^{2k+1}\big]\,\hat\bphi
\\
&=\sum_k\,\hat T^k\hat D\, \big[h_{2k}+h_{2k+1}\,\hat T\big]\,\hat\bphi\,.
\end{array}
\e
The second line in (\ref{tsrb4}) follows from (\ref{tsrb1}), (\ref{qtpd}) and (\ref{pdqt}).
Finally, that the function $\psi$ given by (\ref{tsrb4w}) is an orthonormal wavelet associated with the MRA $\big\{\V_n\big\}_{n\in\Z}$ is a straightforward consequence of Remark \ref{rqmf} with $m=0$.
\end{proof}

\begin{remark}\rm
Not every pair $(h,g)$ of functions of $L^2(\partial {\D})$ verifying conditions (a)--(d) of Proposition \ref{p8} or conditions (1) and (2) of Proposition \ref{p9} comes from an MRA. Under additional conditions,  Cohen \cite{C90} and Lawton \cite{L91} gave equivalent conditions for it.
This question together a deep analysis of relations (\ref{tsrb4}) and (\ref{tsrb4w}) shall be considered elsewhere.
\end{remark}

%---------------

\section*{Acknowledgments}
The authors wish to thank  Ulrich Sorger and Tomasz Ignatz for useful discussions. This work was partially supported by research projects VA108A08 and GR224 (Castilla y Le\'on).
%-----------------------------

\appendix

\section{Some bases for the spectral models}\label{appA}

In order to obtain the respective spectral models given in Propositions \ref{pdft} and \ref{ptft} for the dilation and translation operators, $D$ and $T$, one must consider orthonormal bases (ONB) of $L^2(\R)$,
$$
\{K_{\pm,j}^{(m)}:j\in \J,\,m\in\Z\}\quad \text{ and }\quad \{L_{i}^{(n)}:i\in \I,\,m\in\Z\}\,,
$$
such that $\{K_{\pm,j}^{(0)}:j\in \J\}$ is an ONB of $L^2[\pm 1,\pm 2)$, $\{L_{i}^{(0)}:i\in \I\}$ is an ONB of $L^2[0,1)$ and
$$
K_{\pm,j}^{(m)}=D^m\,K_{\pm,j}^{(0)}\,,\quad 
L_{i}^{(n)}= T^n\,L_{i}^{(0)}\,,\quad (m,n\in\Z)\,.
$$
The change of representation between the two models is governed by the matrix $\big(\alpha_{i,n}^{s,j,m}\big)$ whose elements satisfy
$$
L_i^{(n)}=\sum_{s,j,m} \alpha_{i,n}^{s,j,m}\, K_{s,j}^{(m)}\,.
$$

%---------

\subsection{Exponential bases}\label{appA1}

For the {\it exponential bases}
$$
K_{\pm,j}^{(0)}(x):=e^{2\pi i jx} \quad \text{ and }\quad L_{k}^{(0)}(x):=e^{2\pi ikx}
$$
one has $\I=\J=\Z$ and the change of representation matrix $\big(\alpha_{k,n}^{s,j,m}\big)$ is given by:

$$
\alpha_{k,0}^{s,j,m}=\left\{
\begin{array}{ll}
\ds -\frac{i\,2^{m/2}\,e^{2\pi ik2^{-m}}}{2\pi(k-j2^m)}\,\big(e^{2\pi ik2^{-m}}-1\big),& {\rm if\,\,} s=+,\,\, m>0,\,\,k-j2^m\neq 0,\\
\ds 2^{-m/2},& {\rm if\,\,} s=+,\,\, m>0,\,\,k-j2^m=0,\\
\ds 0, &  {\rm otherwise.}
\end{array}
\right.
$$

$$
\alpha_{k,-1}^{s,j,m}=\left\{
\begin{array}{ll}
\ds \frac{i\,2^{m/2}\,e^{-2\pi ik2^{-m}}}{2\pi(k-j2^m)}\,\big(e^{-2\pi ik2^{-m}}-1\big),& {\rm if\,\,} s=-,\,\, m>0,\,\,k-j2^m\neq 0,\\
\ds 2^{-m/2},& {\rm if\,\,} s=-,\,\, m>0,\,\,k-j2^m=0,\\
\ds 0, &  {\rm otherwise.}
\end{array}
\right.
$$

$$
\alpha_{k,1}^{s,j,m}=\left\{
\begin{array}{ll}
\ds 1,& {\rm if\,\,} s=+,\,\, m=0,\,\,k=j,\\
\ds 0, &  {\rm otherwise.}
\end{array}
\right.
$$

$$
\alpha_{k,-2}^{s,j,m}=\left\{
\begin{array}{ll}
\ds 1,& {\rm if\,\,} s=-,\,\, m=0,\,\,k=j,\\
\ds 0, &  {\rm otherwise.}
\end{array}
\right.
$$

\noindent
For $p=1,2,\ldots$ and $q=0,1,\ldots,2^p-1$,

$$
\alpha_{k,2^p+q}^{s,j,m}=\left\{
\begin{array}{ll}
\ds -\frac{i\,2^{p/2}\,e^{-2\pi ij2^{-p}q}}{2\pi(k2^p-j)}\,\big(e^{-2\pi ij2^{-p}}-1\big),& {\rm if\,\,} s=+,\,\, m=-p,\,\,k2^p-j\neq 0,\\
\ds 2^{-p/2},& {\rm if\,\,} s=+,\,\, m=-p,\,\,k2^p-j= 0,\\
\ds 0, &  {\rm otherwise.}
\end{array}
\right.
$$

$$
\alpha_{k,-2^p-q-1}^{s,j,m}=\left\{
\begin{array}{ll}
\ds \frac{i\,2^{p/2}\,e^{2\pi ij2^{-p}q}}{2\pi(k2^p-j)}\,\big(e^{2\pi ij2^{-p}}-1\big),& {\rm if\,\,} s=-,\,\, m=-p,\,\,k2^p-j\neq 0,\\
\ds 2^{-p/2},& {\rm if\,\,} s=-,\,\, m=-p,\,\,k2^p-j= 0,\\
\ds 0, &  {\rm otherwise.}
\end{array}
\right.
$$

%---------

\subsection{Haar bases}\label{appA2}

Let $\psi$ be the {\it Haar wavelet} and $\varphi$ the corresponding scaling function,
$$
\psi=\chi_{[0,1/2)}-\chi_{[1/2,1)}\,,\quad \varphi=\chi_{[0,1)}\,,
$$
where the characteristic function $\chi$ is defined by
$$
\chi_{[a,b)}(x)=\left\{\begin{array}{ll}1,& \text{ if }x\in[a,b),\\ 0,& \text{ otherwise.}\end{array}\right.
$$
As usual, for each $j,k\in\Z$, let $\psi_{j,k}:=D^jT^k\psi$ and $\varphi_{j,k}:=D^jT^k\varphi$.

Consider the {\it Haar bases} defined by
$$
\begin{array}{l}
K_{+,0}^{(0)}:=\varphi_{0,1}\,;\\
K_{+,2^p+q}^{(0)}:=\psi_{p,2^p+q}\,, \quad \text {for }p=0,1,\,\dots;\,q=0,1,\ldots,2^{p}-1\,;
\\[1ex]
K_{-,0}^{(0)}:=\varphi_{0,-2}\,;\\
K_{-,2^p+q}^{(0)}:=\psi_{p,-2^{p+1}+q}\,, \quad \text {for }p=0,1,\,\dots;\,q=0,1,\ldots,2^{p}-1\,;
\\[1ex]
L_{0}^{(0)}:=\varphi_{0,0}\,;\\
L_{2^p+q}^{(0)}:=\psi_{p,q}\,, \quad \text {for }p=0,1,\,\dots;\, q=0,1,\ldots,2^{p}-1\,.
\end{array}
$$
Then $\I=\J=\N\cup\{0\}$ and the change of representation matrix $\big(\alpha_{i,n}^{s,j,m}\big)$ is as follows:

\smallskip\noindent
For $n=0$, $i=0$, 
$$
\alpha_{0,0}^{s,j,m}=\left\{
\begin{array}{ll}
\ds 2^{-m/2},& \text{ if } s=+,\,\,j=0,\,\, m>0,\\
\ds 0, &  \text{ otherwise.}
\end{array}
\right.
$$
For $n=0$, $r=0,1,2,\ldots$,
$$
\alpha_{2^r,0}^{s,j,m}=\left\{
\begin{array}{ll}
\ds -2^{-1/2},& \text{ if } s=+,\,\, j=0,\,\, m=r+1,\\
\ds 2^{(r-m)/2},& \text{ if } s=+,\,\, j=0,\,\, m>r+1,\\
\ds 0, &  \text{ otherwise.}
\end{array}
\right.
$$
For $n=0$, $r=0,1,2,\ldots$ and $t=2^p+q$ (with $0\leq p<r$ and $q=0,1,\ldots,2^p-1$),
$$
\alpha_{2^r+t,0}^{s,j,m}=\left\{
\begin{array}{ll}
\ds 1,& \text{ if } s=+,\,\, j=t,\,\, m=r-p,\\
\ds 0, &  \text{ otherwise.}
\end{array}
\right.
$$
For $n=1$ and every $i\geq0$,
$$
\alpha_{i,1}^{s,j,m}=\left\{
\begin{array}{ll}
\ds 1,& \text{ if } s=+,\,\, j=i,\,\, m=0,\\
\ds 0, &  \text{ otherwise.}
\end{array}
\right.
$$
For $n>1$ and $i>0$, with $n=2^u+v$ ($u=1,2,\ldots$, $v=0,1,\ldots,2^u-1$) and $i=2^r+t$ ($r=0,1,\ldots$, $t=0,1,\ldots,2^r-1$), 
$$
\alpha_{2^r+t,2^u+v}^{s,j,m}=\left\{
\begin{array}{ll}
\ds 1,& \text{ if } s=+,\,\, j=2^r(2^u+v)+t,\,\, m=-u,\\
\ds 0, &  \text{ otherwise.}
\end{array}
\right.
$$
For $n>1$ and $i=0$, with $n=2^u+v$ ($u=1,2,\ldots$, $v=0,1,\ldots,2^u-1$), 
$$
\alpha_{0,2^u+v}^{s,j,m}=\left\{
\begin{array}{ll}
\ds 2^{-u/2},& \text{ if } s=+,\,\, j=0,\,\, m=-u,\\
\ds (-1)^{w(u,v,p)}\, 2^{(p-u)/2},& \text{ if } s=+,\,\, j=2^p+\left[{v}/{2^{u-p}}\right]\\
&\text{ for }0\leq p<u,\,\, m=-u,\\
\ds 0, &  \text{ otherwise,}
\end{array}
\right.
$$
where $[\cdot]$ denotes ``entire part of" and\footnote{Taking into account the binary expression $v=\sum_{k=0}^{u-1} t_k\,2^k$, with $t_k=0$ or $1$, one has $w(u,v,p)=t_{u-p-1}$ and $\left[{v}/{2^{u-p}}\right]=\sum_{k=u-p}^{u-1} t_k\,2^{k-(u-p)}$.}
\b\label{apph1}
w(u,v,p)=\left[\frac{v-2^{u-p}[v/2^{u-p}]}{2^{u-p-1}}\right].
\e
For $n=-1$, $i=0$, 
$$
\alpha_{0,-1}^{s,j,m}=\left\{
\begin{array}{ll}
\ds 2^{-m/2},& \text{ if } s=-,\,\,j=0,\,\, m>0,\\
\ds 0, &  \text{ otherwise.}
\end{array}
\right.
$$
For $n=-1$, $r=0,1,2,\ldots$,
$$
\alpha_{2^{r+1}-1,-1}^{s,j,m}=\left\{
\begin{array}{ll}
\ds 2^{-1/2},& \text{ if } s=-,\,\, j=0,\,\, m=r+1,\\
\ds -2^{(r-m)/2},& \text{ if } s=-,\,\, j=0,\,\, m>r+1,\\
\ds 0, &  \text{ otherwise.}
\end{array}
\right.
$$
For $n=-1$, $r=1,2,\ldots$, $0\leq p<r$ and $q=0,1,\ldots,2^p-1$,
$$
\alpha_{2^{r+1}-2^{p+1}+q,-1}^{s,j,m}=\left\{
\begin{array}{ll}
\ds 1,& \text{ if } s=-,\,\, j=2^p+q,\,\, m=r-p,\\
\ds 0, &  \text{ otherwise.}
\end{array}
\right.
$$
For $n=-2$ and every $i\geq0$,
$$
\alpha_{i,-2}^{s,j,m}=\left\{
\begin{array}{ll}
\ds 1,& \text{ if } s=-,\,\, j=i,\,\, m=0,\\
\ds 0, &  \text{ otherwise.}
\end{array}
\right.
$$
For $n<-2$ and $i>0$, with $n=-2^{u+1}+v$ ($u=1,2,\ldots$, $v=0,1,\ldots,2^u-1$) and $i=2^r+t$ ($r=0,1,\ldots$, $t=0,1,\ldots,2^r-1$), 
$$
\alpha_{2^r+t,-2^{u+1}+v}^{s,j,m}=\left\{
\begin{array}{ll}
\ds 1,& \text{ if } s=-,\,\, j=2^r(2^u+v)+t,\,\, m=-u,\\
\ds 0, &  \text{ otherwise.}
\end{array}
\right.
$$
For $n<-2$ and $i=0$, with $n=-2^{u+1}+v$ ($u=1,2,\ldots$, $v=0,1,\ldots,2^u-1$),
$$
\alpha_{0,-2^{u+1}+v}^{s,j,m}=\left\{
\begin{array}{ll}
\ds 2^{-u/2},& \text{ if } s=-,\,\, j=0,\,\, m=-u,\\
\ds (-1)^{w(u,v,p)}\, 2^{(p-u)/2},& \text{ if } s=-,\,\, j=2^p+\left[{v}/{2^{u-p}}\right]\\
&\text{ for }0\leq p<u,\,\, m=-u,\\
\ds 0, &  \text{ otherwise,}
\end{array}
\right.
$$
where $w(u,v,p)$ is given by (\ref{apph1}).

%-----------------


\begin{thebibliography}{99}

\bibitem{BS} 
M. S. Birman, M. Z. Solomjak, 
{\it Spectral Theory of Self Adjoint Operators in Hilbert Space}, 
Reidel, Dordrecht, 1987.

\bibitem{C90}
A. Cohen,
{\it Ondelettes, analyses multirésolutions et filtres miroirs en quadrature},
Ann. Inst. H. Poincaré Anal. Non Linéaire 7 (1990), no. 5, 439--459.

\bibitem{D92}
I. Daubechies, 
{\it Ten lectures on wavelets}, CBMS-NSF Regional Conference Series in Applied Mathematics, 61. SIAM, Philadelphia, PA, 1992.

\bibitem{HALMOS61}
P. R. Halmos,
{\it Shifts on Hilbert spaces},
{J. Reine Angew. Math.} 208 (1961) 102--112.

\bibitem{HELSON64}
H. Helson,
{\it Lectures on invariant subspaces}, 
Academic Press, New York, 1964.

\bibitem{HW96}
E. Hern\'andez, G. Weiss,
{\it  A first course on wavelets},
Studies in Advanced Mathematics. CRC Press, Boca Raton, FL, 1996.

\bibitem{HP74}
E. Hille, R. S. Phillips,
{\it  Functional analysis and semi-groups} (Third printing of the revised edition of 1957),  American Mathematical Society, Providence, R.I., 1974.

\bibitem{JMR01}
S. Jaffard, Y. Meyer, R. D. Ryan,
{\it Wavelets. Tools for science \& technology}, SIAM, Philadelphia, PA, 2001.

\bibitem{L91}
 W. M. Lawton,
{\it Necessary and sufficient conditions for constructing orthonormal wavelet bases},
J. Math. Phys. 32 (1991), no. 1, 57--61.

\bibitem{J96}
P. E. T. Jorgensen,
{\it Scattering theory for orthogonal wavelets}, in {\it Clifford algebras in analysis and related topics}, 173--198,
Stud. Adv. Math., CRC, Boca Raton, FL, 1996. 

\bibitem{L61}
P. D. Lax, {\it Translation invariant spaces},  Proc. Internat. Sympos. Linear Spaces, pp. 299--306, Jerusalem Academic Press, Jerusalem, 1961 .

\bibitem{LMR97}
A. K. Louis, P. Maa\ss, A. Rieder,
{\it Wavelets. Theory and Applications}, John Wiley \& Sons, Chichester, 1997.

\bibitem{M98}
S. Mallat,
{\it A wavelet tour of signal processing}, Academic Press, Inc., San Diego, CA, 1998.

\bibitem{P04}
J. A. Packer, 
{\it Applications of the work of Stone and von Neumann to wavelets}, in  {\it Operator algebras, quantization, and noncommutative geometry},  253--279, Contemp. Math., 365, Amer. Math. Soc., Providence, RI, 2004.

\bibitem{S63}
T. P. Srinivasan, {\it Simply invariant subspaces},  Bull. Amer. Math. Soc.  69  (1963) 706--709.

\bibitem{VEO99}
S. J. L. van Eijndhoven, P. J. Oonincx, 
{\it Frames, Riesz systems and multiresolution analysis in Hilbert spaces},  Indag. Math. (N.S.)  10  (1999) 369--382.

\bibitem{VN49}
J. von Neumann,
{On rings of operators: reduction theory},
Ann. of Math. {50} (1949) 401--485.

\end{thebibliography}
\end{document}